\documentclass[]{amsart}
\usepackage{amssymb,amsmath,cases}

\theoremstyle{plain}
\newtheorem{Theorem}{Theorem}[section]
\newtheorem{Lemma}[Theorem]{Lemma}
\newtheorem{Proposition}[Theorem]{Proposition}
\newtheorem{Corollary}[Theorem]{Corollary}
\newtheorem{Remark}[Theorem]{Remark}

\def\f{\varphi}
\def\Proof{\noindent{{\it Proof. }}}
\def\square{\vbox{
	\hrule height .4pt
	\hbox{\vrule width .4pt height 7pt \kern 7pt
	   \vrule width .4pt}
	\hrule height .4pt }}
\def\QED{\hfill {$\square$}\goodbreak \medskip}

\def\QED{\hfill {$\square$}\goodbreak \medskip}

\def\irn{{\int_{\mathbb{R}^n}}}
\def\R{{\mathbb R}}
\def\f{{\varphi}}

\def\eps{{\varepsilon}}

\def\womega{{\Omega'}}

\def\R{\mathbb R}

\def\weak{{\,\rightharpoonup\,}}

\def\pstarstar{{p^{*\!*}}}

\begin{document}

\title{OPTIMAL RELLICH-SOBOLEV CONSTANTS AND THEIR EXTREMALS}

	\author{
	{\sc Roberta Musina}}
	
	\date{}

	\maketitle
\small
	
	\centerline{Dipartimento di Matematica ed Informatica}
	
	\centerline
	{Universit\`a di Udine via delle Scien\-ze, 206 -- 33100 Udine, Italy}
	
	\centerline{
	Email:{musina@dimi.uniud.it}}


\bigskip

\small
\noindent{\bf Abstract.}
		We prove that extremals for second order Rellich-Sobolev inequalities 
have constant sign. Then we show that the optimal constants in Rellich-Sobolev inequalities
on a bounded domain $\Omega$ and under Navier boundary conditions do not depend on $\Omega$.
\normalsize

\section{Introduction and main results}
\let\thefootnote\relax\footnote
{AMS Subject Classifications: 46E35, 26D10, 35J55.}
\let\thefootnote\relax\footnote
{R. Musina is partially supported by Miur-PRIN project 2009WRJ3W7-001 ``Fenomeni di concentrazione e {pro\-ble\-mi} di analisi geometrica".}
Let $n\ge 3$ be an integer and $p,q$ given exponents, such that 
$$
1<p<q~,\quad\text{and ~
$q\le \pstarstar=\displaystyle\frac{np}{n-2p}$ ~if $n>2p$.}
$$
 Assume that  $\alpha,\beta\in\R$ are given in such a way that
\begin{subnumcases} {\label{eq:ma} }
~\beta=n-q\frac{n-2p+\alpha}{p}&{}\label{eq:ma_a}\\
~2p-n<\alpha<np-n.&{}\label{eq:ma_b}
\end{subnumcases}
In \cite{MS} it has been proved that there exists a 
best constant $c>0$ such that
\begin{equation}
\label{eq:inequality}
\irn |x|^{\alpha}|\Delta u|^p dx\ge c\left(\irn |x|^{-\beta}|u|^{q}~\!dx\right)^{p/q}
\end{equation}
for any $u\in C^2_c(\R^n\setminus\{0\})$, see also Corollary \ref{C:deltanabla}
in Subsection \ref{A:D2}.
A rescaling argument plainly shows that (\ref{eq:ma_a}) is a necessary condition for the 
validity of (\ref{eq:inequality}), while  
assumption (\ref{eq:ma_b}) can be weakened. For instance, if
$p=2$ then there exists 
an increasing sequence of integers $j_k\ge n-2$ such that $j_k\to\infty$
and such that
(\ref{eq:inequality}) holds with a  constant $c>0$ if and only if
$(\alpha-2)^2\neq j_k$ for any $k\ge 1$, compare with 
\cite[Theorem 1.1]{CM2}.

When  $q=p>1$, $\beta=2p-\alpha$ and (\ref{eq:ma_b}) hold, then (\ref{eq:inequality}) includes the  
sharp Rellich-type inequality
\begin{equation}
\label{eq:Rellich_Mit_Rn}
\irn |x|^\alpha|\Delta u|^{p}~\!dx\ge 
\gamma_{p,\alpha}^p\irn|x|^{\alpha-2p}|u|^{p}~\!dx \quad\forall u\in C^2_c(\R^n),
\end{equation}
that has been proved by Mitidieri in \cite{Mit2}. Here we have put
\begin{equation}
\label{eq:gamma_intro}
\gamma_{p,\alpha}:=
\frac{n-2p+\alpha}{p}~\!\frac{np-n-\alpha}{p}.
\end{equation}
In particular, inequality
(\ref{eq:inequality}) is naturally related to the function space
$$
\mathcal D^{2,p}(\R^n;|x|^\alpha dx):=
\left\{u\in L^p(\R^n;|x|^{\alpha-2p}dx)~|~\Delta u\in L^p(\R^n;|x|^\alpha dx)~\right\}~\!,
$$
and one is lead to study the
minimization problem
\begin{equation}
\label{eq:space_weight}
S_{p,q}(\alpha)=\inf_{\scriptstyle u\in \mathcal D^{2,p}(\R^n;|x|^\alpha dx)\atop\scriptstyle  u\ne 0}
\frac{\displaystyle\irn |x|^{\alpha}|\Delta u|^p dx}
{\left(\displaystyle\irn |x|^{-\beta}|u|^{q}~\!dx\right)^{p/q}}~.
\end{equation}

If  $n>2p$, $q=\pstarstar$ and $\alpha=\beta=0$, then 
the infimum in (\ref{eq:space_weight}) equals
 the Sobolev constant
\begin{equation}
\label{eq:Sob}
S_p=\inf_{\scriptstyle u\in \mathcal D^{2,p}(\R^n)\atop\scriptstyle u\ne 0}
\frac{\displaystyle\irn|\Delta u|^p~\!dx}
{\left(\displaystyle\irn |u|^{\pstarstar}~\!dx\right)^{p/\pstarstar}}~,
\end{equation}
that is relative to the critical embedding $\mathcal D^{2,p}(\R^n)\hookrightarrow L^{\pstarstar}(\R^n)$.
In \cite{PLL}, Corollary I.2, P.L. Lions proved  that every bounded minimizing sequence is
relatively compact up to dilations and translations, and in particular $S_p$ is achieved.
Moreover,
by using Schwarz symmetrization he showed that,  up to a change of sign,
any extremal for $S_p$ is spherically symmetric, positive and decreasing.
Using this information,  Hulshof and Van der Vorst \cite{HVdV} were able to prove  uniqueness
of  extremals for $S_p$,
modulo dilations, translations in $\R^n$ and change of sign.

As concerns general exponents $p,q,\alpha$ and $\beta$ satisfying (\ref{eq:ma}),
one can find sufficient conditions for the existence of extremals for
$S_{p,q}(\alpha)$
in the appendix of \cite{MS}, see also \cite{CM2} for the Hilbertian case $p=2$.

In presence of weights rearrangement techniques are in general not 
applicable. As a matter of fact, {\em breaking symmetry} may occur, see
Section 5 of \cite{CM2}, where $p=2$ is assumed.
Actually {\em breaking positivity} phenomena can not be a priori excluded as well: indeed,
it may happen that no extremal
for $S_{p,q}(\alpha)$ has constant sign. 
From a merely technical point of view,
the reason lies in the fact that truncation $u\mapsto u^+$ can not be used
to prove the positivity of extremals for $S_{p,q}(\alpha)$, since 
in general $u^+\notin  \mathcal D^{2,p}(\R^n;|x|^\alpha dx)$ for
$u\in  \mathcal D^{2,p}(\R^n;|x|^\alpha dx)$. Nevertheless, in Section \ref{SS:weights} we
prove the next result, that gives a  positive answer to a query raised in 
\cite{CM2}.

\begin{Theorem}
\label{T:noBP}
Assume that (\ref{eq:ma}) holds.
Let $u\in \mathcal D^{2,p}(\R^n;|x|^\alpha dx)$
be an extremal for $S_{p,q}(\alpha)$. Then (up to a change of sign), 
$u$ is positive and superharmonic.
\end{Theorem}

Assumption 
(\ref{eq:ma_b}) can not be neglected:  by
Theorem 4.1 in \cite{CM2}, breaking positivity
does occur if $p=2$, $|\alpha-2|$ is large enough and $q\approx 2$.

\medskip

Next, let $\Omega$ be a bounded and smooth domain containing the origin. Let 
$$
\mathcal D^{2,p}_N(\Omega;|x|^\alpha dx)=\left\{u\in L^p(\Omega;|x|^{\alpha-2p}dx)~|~
\Delta u\in L^{p}(\Omega;|x|^\alpha dx)~,~u= 0~~\text{on $\partial\Omega$}~\right\}
$$
(see Section \ref{SS:bounded} for details). The
optimal Rellich-Sobolev constant  under Navier boundary conditions is given by
\begin{equation}
\label{eq:inf_omegaa}
S_{p,q}^{\rm Nav}(\Omega;\alpha)=\inf_{\scriptstyle u\in \mathcal D^{2,p}_N(\Omega;|x|^\alpha dx)\atop\scriptstyle  u\ne 0}
\frac{\displaystyle\int_\Omega |x|^{\alpha}|\Delta u|^p dx}
{\left(\displaystyle\int_\Omega |x|^{-\beta}|u|^{q}~\!dx\right)^{p/q}}~\!.
\end{equation}
A rescaling argument plainly shows
that $S_{p,q}^{\rm Nav}(\Omega;\alpha)\le S_{p,q}(\alpha)$. The opposite inequality
is not trivial at all, as in general a function $u\in \mathcal D^{2,p}_N(\Omega;|x|^\alpha dx)$ 
can not  be extended to $\overline u\in \mathcal D^{2,p}(\R^n;|x|^\alpha dx)$
by putting $\overline u\equiv 0$ outside $\Omega$.
In Section \ref{SS:alpha_Navier} we prove
the next result.

\begin{Theorem}
\label{T:RS_Omega} If (\ref{eq:ma}) holds, and if $\Omega$ is a bounded domain of
class $C^2$ containing the origin, then
\begin{equation}
\label{eq:weight_Navier}
S_{p,q}^{\rm Nav}(\Omega;\alpha)=S_{p,q}(\alpha)
\end{equation}
and in particular $S_{p,q}^{\rm Nav}(\Omega;\alpha)$ is positive.
Moreover, 
$S_{p,q}^{\rm Nav}(\Omega;\alpha)$ is not achieved in $\mathcal D^{2,p}_{N}(\Omega;|x|^\alpha dx)$.
\end{Theorem}
If $n>2p$, $q=\pstarstar$ and $\alpha=\beta=0$, then
the infimum in (\ref{eq:inf_omegaa}) coincides with 
$$
S_p^{\rm Nav}(\Omega):=\inf_{\scriptstyle u\in W^{2,p}\cap W^{1,p}_0(\Omega)\atop\scriptstyle u\ne 0}
\frac{\displaystyle\int_\Omega|\Delta u|^p~\!dx}
{\left(\displaystyle\int_\Omega |u|^{\pstarstar}~\!dx\right)^{p/\pstarstar}}.
$$
Hence, by Theorem \ref{T:RS_Omega} we have that
\begin{equation}
\label{eq:equality} 
S_p^{\rm Nav}(\Omega)= S_p
\end{equation}
and $S_p^{\rm Nav}(\Omega)$ is not achieved.
In the Hilbertian case $p=2$,  equality (\ref{eq:equality})
has been proved by Van der Vorst in
\cite{VdV} and by Ge in \cite{Ge}.  The general case
$p>1$ has been recently exploited by Gazzola, Grunau and Sweers in \cite{GGS}. 
All the above mentioned papers are based again on a rearrangement argument
that, in general, fails in presence of weights. Our arguments to check the more general equality
(\ref{eq:weight_Navier}) are simpler and self-contained. 

\medskip

In the last theorem we provide an unexpected result.
We denote by $C^2_N(\overline\Omega\setminus\{0\})$  the set of functions $u\in C^2(\overline\Omega)$
such that $u=0$ on $\partial\Omega$ and in a neighborhood of $0$.

\begin{Theorem}
\label{T:Rellich}
Let $\Omega$ be a bounded  domain of class $C^2$ containing the origin. 
Let $q\ge p>1$, $\alpha\in\R$ and define $\beta$ as in (\ref{eq:ma_a}). If $\alpha\ge np-n$, then $$
S_{p,q}^{\rm Nav}(\Omega;\alpha):= \inf_{\scriptstyle u\in  C^2_N(\overline\Omega\setminus\{0\})\atop\scriptstyle  u\ne 0}\frac{\displaystyle
\int_\Omega |x|^{\alpha}|\Delta u|^p dx}
{\left(\displaystyle\int_\Omega |x|^{-\beta}|u|^{q}dx\right)^{p/q}}
=0.
$$
\end{Theorem}

To comment Theorem \ref{T:Rellich} we define also
\begin{gather*}
S_{p,q}^{\rm Dir}(\alpha;\Omega):= \inf_{\scriptstyle u\in  C^2_c(\Omega\setminus\{0\})\atop\scriptstyle  u\ne 0}\frac{\displaystyle
\int_\Omega |x|^{\alpha}|\Delta u|^p dx}
{\left(\displaystyle\int_\Omega |x|^{-\beta}|u|^{q}dx\right)^{p/q}}\\
S_{p,q}(\alpha):=\inf_{\scriptstyle u\in C^2_c(\R^n\setminus\{0\})\atop\scriptstyle  u\ne 0}\frac{\displaystyle
\irn |x|^{\alpha}|\Delta u|^p dx}
{\left(\displaystyle\irn |x|^{-\beta}|u|^{q }dx\right)^{p/q}}.
\end{gather*}
Thanks to Theorems \ref{T:RS_Omega} and \ref{T:Rellich} (see also Section \ref{A:D2}),
we have that
\begin{eqnarray}
\nonumber
S_{p,q}^{\rm Nav}(\Omega;\alpha)= S_{p,q}^{\rm Dir}(\alpha;\Omega)=
S_{p,q}(\alpha)&&\quad\text{if $\alpha\in (2p-n,np-n)$}\\ 
0=S_{p,q}^{\rm Nav}(\Omega;\alpha)\le S_{p,q}^{\rm Dir}(\alpha;\Omega)=
~S_{p,q}(\alpha)&&\quad\text{if $\alpha\ge np-n$.}\label{eq:diff_b}
\end{eqnarray}
In general, the strict inequality holds in (\ref{eq:diff_b}). 
Assume $q\le p^{*\!*}$ if $n>2p$.
By Corollary
\ref{C:deltanabla} below, the infimum
$S_{p,q}(\alpha)$ is positive if and only if $-\gamma_{p,\alpha}$ is not an eigenvalue
of the Laplace-Beltrami operator
on the sphere. Therefore, if in addition it holds that 
$\alpha> np-n$, then
$$
0=S_{p,q}^{\rm Nav}(\Omega;\alpha)<S_{p,q}^{\rm Dir}(\Omega;\alpha)=S_{p,q}(\alpha).
$$
In particular the optimal Rellich-Sobolev constant under Navier boundary conditions never
depends on the domain, but it may not
coincide with the Rellich-Sobolev constant under Dirichlet boundary conditions,
nor with the Rellich-Sobolev constant on the whole space.

The picture is not complete if $\alpha< 2p-n$; see Proposition
 \ref{P:CM1} for a partial result.

\medskip

Our proofs are based on some knowledge of the 
weighted Sobolev space involved, compare with Section \ref{S:spaces},
and on  two basic
facts. The first one concerns  inequality
(\ref{eq:Rellich_Mit_Rn}) and its generalization to Rellich type inequalities on bounded
domains with Navier boundary conditions.
Actually Lemma \ref{L:Mit} in Section \ref{A:D2}  has been already proved in 
\cite{Mit2}, but it has never been explicitly stated in the form we need
for our purposes. 
The second basic fact is the core of Theorem \ref{T:linear} in Section \ref{S:linear}, that is concerned
with non-homogeneous equations
of the form
$$
-\Delta v=f\quad\text{in $\R^n$,}
$$
for $f$ varying in weighted $L^p$ spaces. 
The proofs of the main Theorems \ref{T:noBP}, \ref{T:RS_Omega}
and \ref{T:Rellich}  can be found in the last three sections.

\section{Notation and weighted Sobolev spaces}
\label{S:spaces}

The characteristic function of a domain $\Omega$ in $\R^n$ is denoted by $\chi_\Omega$.
If $f:\Omega\to \R$  is given, then $\chi_\Omega f$ denotes the extension of $f$
by the null function outside $\Omega$. 

We will use several function spaces on $\Omega$. In addition to the usual spaces
$C^k_c(\Omega)$ and $C^\infty_c(\Omega)$, let us denote
\begin{gather*}
C^2_N(\overline\Omega)=\left\{u\in C^2(\overline\Omega)~|~\text{supp$(u)$ is compact, and $u\equiv 0$ on $\partial\Omega$}\right\}\\
C^2_N(\overline\Omega\setminus\{0\})=\left\{u\in C^2_N(\overline\Omega)~|~u\equiv 0~~\text{in a neighborhood of $0$}~\right\}.
\end{gather*}

Let $a\in\R$.  
The weighted Lebesgue space $L^q(\Omega;|x|^a~\! dx)$ is the space of measurable 
$u$ on $\Omega$ having finite norm $\left(\int_{\Omega}|x|^a|u|^q~\!dx\right)^{1/q}$.  For $a=0$ we  write
$L^q(\Omega)$, as usual.

\medskip

A  proper function $u:\Omega\to \R\cup\{\infty\}$ is {\em superharmonic} if it is lower semicontinuous
and 
$$
u(x)\ge \frac{1}{|\partial B_r|}\int_{\partial B_r(x)}u(y)d\sigma(y)
\quad\text{for any $x\in\Omega$, $r<\text{dist}(x,\partial\Omega)$.}
$$
Every superharmonic function
in $\Omega$ belongs to $L^1_{\rm loc}(\Omega)$. Moreover, $u\in L^1_{\rm loc}(\Omega)$ is superharmonic if
and only if $-\Delta u\ge 0$ in  $\mathcal D'(\Omega)$,
that is,
$$
\int_\Omega u(-\Delta\f)~\!dx\ge 0\quad\text{for any nonnegative $\f\in C^\infty_c(\Omega)$.}
$$

Let 
$\Omega$ be a bounded domain with $\partial\Omega$ of class $C^2$. 
By definition, $W^{k,p}_0(\Omega)$ is the closure of $C^k_c(\Omega)$ in the standard Sobolev space $W^{k,p}(\Omega)$. Hence,
 $W^{k,p}_0(\Omega)$ is isometrically embedded into $W^{k,p}(\R^n)$ via the
 null extension $u\mapsto \chi_\Omega u$.

We adopt the  notation
$$
W^{2,p}_{N}(\Omega)=W^{1,p}_{0}(\Omega)\!\cap\!W^{2,p}(\Omega).
$$
It turns out that $W^{2,p}_{N}(\Omega)$ is a Banach space with respect to the norm 
$$
\|u\|=\left(\int_\Omega|\Delta u|^{p}~\!dx\right)^{1/p},
$$
which is equivalent to norm induced by $W^{2,p}(\Omega)$. For $p=2$ we
will use the simplified notation $H^k$ instead of $W^{k,2}$. 

Any function
$u\in W^{2,p}_{N}(\Omega)$ is the limit in $W^{2,p}(\Omega)$ of a
sequence $u_{h}\in C^2_N(\overline\Omega)$. For instance, one can
define $u_h$ to be the projection of 
$(\chi_\Omega u)*\rho_{\varepsilon_h}\in H^1(\R^n)$ on $H^1_0(\Omega)$,
where $\varepsilon_h\to 0^+$ and $\rho_{\varepsilon}$ is the standard $\varepsilon$-mollifier. 

For $n>2p$ let $\mathcal D^{2,p}(\R^{n})$ be the completion of $C^{\infty}_{c}(\R^{n})$ with respect to the norm 
$$\|u\|^p=\irn|\Delta u|^p~\!dx. 
$$
It is well known that the space $\mathcal D^{2,p}(\R^n)$ is continuously embedded into
$L^\pstarstar(\R^n)$. The Sobolev constant
$S_p$  in (\ref{eq:Sob}) is positive and achieved in $\mathcal D^{2,p}(\R^n)$ , see for instance \cite{PLL}, \cite{Sw}.

The next lemma is  
based on a standard trick.

\begin{Lemma}
\label{L:extension1}
Let $\Omega, \womega$ be bounded  domains such that 
$\overline\Omega\subset\womega$. In addition, assume that
$\partial\Omega$
is Lipschitz and
$\partial\womega$ is of class $C^2$. Let 
$u\in W^{2,p}_{N}(\Omega)$. The problem
\begin{equation}
\label{eq:uffa}
\begin{cases}
-\Delta v=\chi_\Omega|-\Delta u|&\text{in $\womega$,}\\
v=0&\text{on $\partial\womega$,}
\end{cases}
\end{equation}
has a unique solution $v \in W^{2,p}_N(\womega)$, and
$$
v\ge  \chi_\Omega|u|
\quad\text{in $\womega$}.
$$
\end{Lemma}

\Proof
We can assume that $u\neq 0$.
The existence of a unique $v\in W^{2,p}_N(\womega)$ solving (\ref{eq:uffa})
is indeed a well known fact. Since $v$ is superharmonic then $v>0$ in $\womega$.
Now $v\pm u\in L^1(\Omega)$ and $-\Delta(v\pm u)\ge 0$ in $\Omega$. Since $v\pm u>0$
on $\partial \Omega$ in the sense of traces, then $v\ge \mp u$ in $\Omega$. Thus $v\ge |u|$ in $\Omega$.
\QED

One of the main tools in our arguments is
the  Hardy inequality in \cite{H}, \cite{HLP}. 
For any $p>1$, $a\in \R$ it holds that
$$
\irn|x|^a|\nabla u|^p~dx\ge|H_{1,a}|^p\irn|x|^{a-p}|u|^p~\!dx
\quad\forall u\in C^1_c(\R^n\setminus\{0\}),
$$
where
\begin{equation}
\label{eq:H1a}
H_{1,a}:=\frac{n+a}{p}-1.
\end{equation}
Moreover, the constant $|H_{1,a}|^p$ can not be improved 
and it is
not achieved in any reasonable function space. 

In the remaining part of this section we describe the weighted Sobolev
spaces that are needed to prove our main results.

\subsection{The spaces $\mathcal D^{1,p}(\R^n;|x|^a dx)$}
\label{A:D1}
~\\
Assume $a\neq p-n$ and define the space $\mathcal D^{1,p}(\R^n;|x|^a dx)$ as the completion of $C^1_c(\R^n\setminus\{0\})$ with respect to the norm
$$
\|u\|=\left(\irn|x|^a|\nabla u|^p~dx\right)^{1/p}.
$$
Then $\mathcal D^{1,p}(\R^n;|x|^a dx)\hookrightarrow L^{p}(\R^n;|x|^{a-p} dx)$ by the Hardy inequality.
\begin{Remark}
\label{R:L1}
Assume $a<pn-n-p$, $a\neq p-n$. Then $\mathcal D^{1,p}(\R^n;|x|^a dx)\subset
L^1_{\rm loc}(\R^n)$.
Indeed, if $u\in \mathcal D^{1,p}(\R^n;|x|^a dx)$ then  
for any bounded domain $\Omega$ we have that
$$
\int_{\Omega}|u|~\!dx\le \left(\int_{\R^n}|x|^{a-p}|u|^p~\!dx\right)^{\frac{1}{p}}
\left(\int_{\Omega}|x|^{\frac{p-a}{p-1}}~\!dx\right)^{\frac{p-1}{p}}<\infty.
$$
\end{Remark}
In order to simplify the proofs it is convenient to 
 introduce the cylinder
$$
 \mathcal Z^n= \R\times\mathbb S^{n-1},
 $$
whose points are denoted by $(s,\sigma)$, and the transform
$$
\mathcal T_{1,a}:C^1_c( \mathcal Z^n)\to C^1_c(\R^n\setminus\{0\})~,\quad
(\mathcal T_{1,a}g)(x)=|x|^{-H_{1,a}}g\left(-\log|x|,\frac{x}{|x|}\right).
$$
The next lemma has been already pointed out in \cite{M}, in a radial setting.

\begin{Lemma}
\label{L:description1}
If $a\neq p-n$, then
$$
\|g\|_{1,a}:=\left(\irn|x|^a|\nabla(\mathcal T_{1,a} g)|^p~dx\right)^{1/p}~\quad\text{for $g\in C^1_c( \mathcal Z^n)$,}
$$
is equivalent to the standard norm in $W^{1,p}(\mathcal Z^n)$. Thus $\mathcal T_{1,a}$
can be uniquely extended to an isomorphism  $W^{1,p}(\mathcal Z^n)\to\mathcal D^{1,p}(\R^n;|x|^a dx)$, and
$$
\mathcal D^{1,p}(\R^n;|x|^a dx)=\left\{u\in L^{p}(\R^n;|x|^{a-p} dx)~|~|\nabla u|\in L^p(\R^n;|x|^{a}dx)~\right\}.
$$
If in addition $a>p-n$,  then $C^1_c(\R^n)\subset \mathcal D^{1,p}(\R^n;|x|^a dx)$.
\end{Lemma}

\Proof
Notice that
\begin{gather*}
\irn|x|^{a-p}|\mathcal T_{1,a} g|^p~\!dx=\int_{\mathcal Z^n}|g|^p~\!dsd\sigma\\
\irn|x|^a|\nabla(\mathcal T_{1,a} g)|^p~\!dx=\int_{\mathcal Z^n}\left|(g_s+H_{1,a}g)^2+|\nabla_\sigma g|^2\right|^{\frac{p}{2}}~\!dsd\sigma
\end{gather*}
for any $g\in C^1_c( \mathcal Z^n)$. In particular, the Hardy inequality and a density argument give
\begin{equation}
\label{eq:cor_Hardy}
\|g\|_{1,a}^p\ge |H_{1,a}|^p\int_{\mathcal Z^n}|g|^p~\!dsd\sigma
\end{equation}
for any $g\in W^{1,p}(\mathcal Z^n)$.
It is easy to prove that 
$\|\cdot\|_{1,a}$ is uniformly bounded from above by the standard norm in $W^{1,p}(\mathcal Z^n)$. To prove the converse
take a sequence $g_h$ such that $\|g_h\|_{1,a}\to 0$. Then $g_h\to 0$ in $L^p(\mathcal Z^n)$ by (\ref{eq:cor_Hardy})
and since $H_{1,a}\neq 0$. Thus
\begin{eqnarray*}
o(1)=\|g_h\|_{1,a}^p&=& \int_{\mathcal Z^n}\left|(g_h)_s^2+|\nabla_\sigma g_h|^2\right|^{\frac{p}{2}}~\!dsd\sigma
+o(1)   \\         
&\ge&
\frac{1}{2}\int_{\mathcal Z^n}\left(|(g_h)_s|^p+|\nabla_\sigma g_h|^p\right)~\!dsd\sigma+o(1),
\end{eqnarray*}
hence $g_h\to 0$ in $W^{1,p}(\mathcal Z^n)$. The equivalence of the two norms is proved. 
To conclude, recall that
$\displaystyle{
W^{1,p}( \mathcal Z^n)=\left\{g\in L^p( \mathcal Z^n)~|~g_s,|\nabla_{\sigma}g|\in
L^p( \mathcal Z^n)~\right\}}$, and notice that the weights $|x|^a, |x|^{a-p}$
are locally integrable if $a>p-n$.
\QED

\begin{Remark}
\label{R:symmetry1}
For $a\neq p-n$, $u\in C^1_c(\R^n\setminus\{0\})$ we put $\hat a=2(p-n)-a$ and
$$
\hat u(x)=u\left(\frac{x}{|x|^2}\right),
$$
respectively.  By direct computation one gets that
$$
\int_{\R^{n}}|x|^{a}|\nabla\! u|^p~\!dx=\int_{\R^{n}}|x|^{\hat a}|\nabla\! \hat u|^p~\!dx.
$$
Thus the functional transform $u\mapsto \hat u$ can be extended to a unique isometry
$$\mathcal D^{1,p}(\R^n;|x|^a dx)\to\mathcal D^{1,p}(\R^n;|x|^{\hat a} dx).$$
\end{Remark}

In the next result we provide an alternative proof of the celebrated
Maz'ya and Caffarelli-Kohn-Nirenberg inequalities in  \cite{Ma}, \cite{CKN}. 

\begin{Lemma}
\label{L:CKN}
Let $q>p$, and assume $q\le p^*=\frac{np}{n-p}$ if $n>p$. If $a\neq p-n$  then there exists a positive best  constant
$s_{p,q}(a)$ such that
$$
\irn|x|^a|\nabla u|^p~\!dx\ge s_{p,q}(a)\left(\irn|x|^{-n+qH_{1,a}}|u|^q~\!dx\right)^{p/q}
$$
for any $u\in \mathcal D^{1,p}(\R^n;|x|^{a}dx)$. In addition, $s_{p,q}(a)=s_{p,q}(2(p-n)-a)$.
\end{Lemma}

\Proof Notice  that
$$
\irn|x|^{-n+qH_{1,a}}|\mathcal T_{1,a}g|^q~\!dx=\int_{\mathcal Z^n}|g|^q~\!dsd\sigma
$$
for any $g\in C^1_c(\mathcal Z^n)$. Since 
$W^{1,p}(\mathcal Z^n)\hookrightarrow L^q(\mathcal Z^n)$ by Sobolev embedding theorem, we readily
infer that $\mathcal D^{1,p}(\R^n;|x|^{a}dx)\hookrightarrow L^q(\R^n;|x|^{-n+qH_{1,a}}dx)$
with a continuous embedding,
and the desired inequality follows. Finally, $s_{p,q}(a)=s_{p,q}(2(p-n)-a)$ by
Remark \ref{R:symmetry1}.
\QED

\begin{Remark}
Assume $n>p$ and take $a=0$. The above lemmata apply to the
standard space $ \mathcal D^{1,p}(\R^n)$. In particular, $C^1_c(\R^n\setminus\{0\})$ is dense in
$ \mathcal D^{1,p}(\R^n)$, and $ \mathcal D^{1,p}(\R^n)$ can be identified with
$W^{1,p}(\mathcal Z^n)$
via the transform $\mathcal T_{1,a}$. These facts are well known when $p=2$, see for instance
\cite{CaWa}.
\end{Remark}

The next  maximum principle for superharmonic functions
 might have an independent interest. 
 
\begin{Theorem}
\label{T:MM}
Assume $p-n<a< pn-n-p$. If $\omega\in \mathcal D^{1,p}(\R^n;|x|^{a}dx)\setminus\{0\}$ 
is superharmonic in $\R^n\setminus\{0\}$, then $\omega$
is superharmonic and strictly positive on $\R^n$.
\end{Theorem}

\Proof
First of all we recall that $\omega\in L^1_{\rm loc}(\R^n)$ by Remark \ref{R:L1}.
Fix any nonnegative $\eta\in C^\infty_c(\R^n)$.
Notice that $-\frac{a}{p-1}>p'-n$, where $p'=\frac{p}{p-1}$. Use Lemma \ref{L:description1}
to infer that 
 $\eta\in \mathcal D^{1,p'}(\R^n;|x|^{-\frac{a}{p-1}}dx)$. Thus there exists a sequence 
 $\eta_h\in C^\infty_c(\R^n\setminus\{0\})$ such that 
$\eta_h\to \eta$ in $\mathcal D^{1,p'}(\R^n;|x|^{-\frac{a}{p-1}}dx)$.  Using
truncation and a standard  convolution argument we can assume that $\eta_h\ge 0$. 
Since $\omega\in W^{1,p}_{\rm loc}(\R^n\setminus\{0\})$ is superharmonic on $\R^n\setminus\{0\}$, then
$$
0\le\irn\omega(-\Delta \eta_h)~\!dx=\irn\nabla\omega\cdot\nabla\eta_h~\!dx=\irn\left(|x|^{\frac{a}{p}}\nabla \omega\right)
\cdot \left(|x|^{-\frac{a}{p}}\nabla \eta_h\right)~\!dx.
$$
Next notice that
$$
|x|^{\frac{a}{p}}\nabla \omega\in L^p(\R^n)^n~,\quad
|x|^{-\frac{a}{p}}\nabla \eta_h\to |x|^{-\frac{a}{p}}\nabla \eta\quad
\text{in $L^{p'}(\R^n)^n$}.
$$
Therefore we can pass to the limit to infer
$$
0\le\irn\nabla\omega\cdot\nabla\eta~\!dx=\irn\omega(-\Delta\eta)~\!dx.
$$ 
Thus $-\Delta\omega\ge 0$ as a distribution in $\mathcal D'(\R^n)$, as $\eta$ was
arbitrarily chosen.

To conclude the proof we only have to show that $\omega\ge 0$ almost everywhere in $\R^n$.
For sake of clarity we first assume  $p\ge 2$, as the proof needs less computations in this case. 
Use Lemma \ref{L:description1} and known results on truncations  to get that
 $\omega^-:=-\min\{\omega,0\}\in \mathcal D^{1,p}(\R^n;|x|^{a}dx)$. Then approximate
 $\omega^-$ in $\mathcal D^{1,p}(\R^n;|x|^{a}dx)$ with a sequence of functions in $C^1_c(\R^n\setminus\{0\})$
 to prove that we can test  $-\Delta\omega\ge 0$ with $|x|^{a+2-p}(\omega^-)^{p-1}$.
 At the end one gets
 \begin{eqnarray*}
\irn(-\Delta\omega)|x|^{a+2-p}(\omega^-)^{p-1}~\!dx&=&
\irn\nabla\omega\cdot\nabla\left(|x|^{a+2-p}(\omega^-)^{p-1}\right)~\!dx\\
&=&-\irn\nabla\omega^-\cdot\nabla\left(|x|^{a+2-p}(\omega^-)^{p-1}\right)~\!dx\ge 0,
\end{eqnarray*}
that is,
\begin{eqnarray*}
(p-1)\irn |x|^{2+a-p}|\nabla \omega^-|^2(\omega^-)^{p-2} dx&\le&
-\frac{1}{p}\irn\nabla|x|^{2+a-p}\cdot\nabla(\omega^-)^p dx\\
&=&H_{1,a}(2+a-p)\irn|x|^{a-p}|\omega^-|^p dx.
\end{eqnarray*}
If $a\le p-2$ we readily get that $\omega^-\equiv 0$. Otherwise, we
use the Hardy inequality
$$
\irn|x|^{2+a-p}|\nabla v|^2~\!dx\ge \left|\frac{n+a-p}{2}\right|^2\irn|x|^{a-p}|v|^2~\!dx
$$
with $v=(\omega^-)^{\frac{p}{2}}\in \mathcal D^{1,2}(\R^n;|x|^{2+a-p}dx)$,
to infer
\begin{eqnarray*}
(p-1)\irn |x|^{2+a-p}|\nabla \omega^-|^2|\omega^-|^{p-2} dx&\le&
\frac{4(2+a-p)}{p(n+a-p)}\irn|x|^{2+a-p}\left|\nabla|\omega^-|^{\frac{p}{2}}\right|^2 dx\\
&=& \frac{2+a-p}{H_{1,a}}\irn|x|^{2+a-p}|\omega^-|^{p-2}|\nabla\omega^-|^2 dx
\end{eqnarray*}
Since $H_{1,a}^{-1}(2+a-p)<p-1$ as $a<np-p-n$,
then necessarily $\omega^-\equiv 0$, that is, $\omega\ge 0$. If $p\in(1,2)$ one repeats the same argument,
with $(\omega^-)^{p-2}\omega^-$ replaced by
$$
\omega_{\eps}=\left(|\omega|^2+\eps^2\right)^{\frac{p-2}{2}}\omega^-,
$$
where $\eps\to 0^+$, to get in similar way that $\omega\ge 0$ almost everywhere on $\R^n$.
The proof is complete, as every superharmonic, nontrivial and nonnegative
function on $\R^n$ is everywhere positive.
\QED

\subsection{The spaces $\mathcal D^{2,p}(\R^n;|x|^\alpha dx)$}
\label{A:D2}
~\\
For {\bf any} given exponent $\alpha\in \R$ we introduce the weighted Rellich constant
$$
 \mu_{p,\alpha}:  =  \inf_{\scriptstyle u\in C^{2}_{c}(\R^n\setminus\{0\})\atop\scriptstyle  u\ne 0}\frac{\displaystyle
\irn |x|^{\alpha}|\Delta u|^p dx}
{\displaystyle\irn |x|^{\alpha-2p}|u|^{p }dx}~\!.
$$
A crucial role
is played by the constants
\begin{equation}
\label{eq:H}
H_{2,\alpha}=\frac{n+\alpha}{p}-2~,\quad
\gamma_{p,\alpha}:=
\left(n-\frac{n+\alpha}{p}\right)H_{2,\alpha},
\end{equation}
see also (\ref{eq:gamma_intro}).
In  \cite{MSS}, Metafune, Sobajima and Spina proved that
\begin{equation}
\label{eq:MSS}
\mu_{p,\alpha}>0\quad\text{ if and only if \quad$-\gamma_{p,\alpha}\notin\left\{ k(n-2+k)~
:~ k\in\mathbb N\cup\{0\}~\right\}$},
\end{equation}
solving a problem that has been left open for long time. 

\begin{Remark}
\label{R:resonance} It is easy to check that
$\mu_{p,\alpha}=0$ if  $-\gamma_{p,\alpha}=k(n-2+k)$ for an integer $k\ge 0$.
For the proof, let $\f_k\in H^1(\mathbb S^{n-1})$ be an eigenfunction of $-\Delta_\sigma$
(the Laplace-Beltrami operator
on the sphere) relative to the eigenvalue $\lambda_k=k(n-2+k)$. 
Fix a nontrivial function $\omega\in C^2_c(\R_+)$ and for any $\eps>0$ use polar coordinates
$(r,\sigma)\in\R_+\times\mathbb S^{n-1}$ to define 
$$
u_\eps(r\sigma):=r^{-H}\omega(r^\eps)\f_k(\sigma)~\!,
$$ 
where we have put $H=H_{2,\alpha}$ to simplify notation.
Then  test $\mu_{p,\alpha}$ with $u_\eps$. By direct computation one gets
\begin{eqnarray*}
\mu_{p,\alpha}&\le&\eps^p~\frac
{\displaystyle\int_0^\infty s^{p-1}|\eps s\omega''+(n-2-2H+\eps)\omega'|^p~\!ds}
{\displaystyle\int_0^\infty s^{-1}|\omega|^p~\!ds}.
\end{eqnarray*}
The conclusion follows by taking the limit as $\eps\to 0$.
\end{Remark}
The explicit value of $\mu_{2,\alpha}$ (case $p=2$) has been computed in \cite{GM}, \cite{CM1}
and \cite{MSS}. 
The sharp value of $\mu_{p,\alpha}$ in case of general exponents $\alpha, p$ is not known yet, unless
$\gamma_{p,\alpha}$ is positive (hence,
$-\gamma_{p,\alpha}$ is below the spectrum of $-\Delta_\sigma$). The weighted Rellich inequality in the next Lemma has been essentially proved in \cite{Mit2}. We cite also \cite{GGM}, Lemma 2,
where $p=2$ is assumed.
We sketch  its proof for the convenience of the reader.

\begin{Lemma}
\label{L:Mit}
Let $\Omega$ be a  domain in $\R^n$. If $\partial \Omega$ is not empty, assume that $\partial \Omega$ is of class $C^2$.
Let  $p>1$ and $\alpha\in\R$ such that $2p-n<\alpha<np-n$. Then 
\begin{equation}
\label{eq:Rellich_Mit_Omega}
\gamma_{p,\alpha}^p\int_{\Omega}|x|^{\alpha-2p}|u|^{p}~\!dx\le \int_\Omega |x|^\alpha|\Delta u|^{p}~\!dx\quad\forall u\in C^2_N(\overline\Omega).
\end{equation}
In particular, we have that (\ref{eq:Rellich_Mit_Rn}) holds.
If $0\in\Omega$ then the constant in the left hand side of (\ref{eq:Rellich_Mit_Omega}) can not be improved.
\end{Lemma}

\Proof
If $p\ge 2$ then (\ref{eq:Rellich_Mit_Omega}) is an immediate consequence of H\"older and Hardy inequalities and of the identity
\begin{eqnarray*}
&&\int_\Omega(-\Delta u) |x|^{\alpha-2p+2}|u|^{p-2}u~\!dx\\
&&\quad\quad=\frac{4(p-1)}{p}\int_\Omega|x|^{\alpha-2p+2}\left|\nabla|u|^{\frac{p}{2}}\right|^2~\!dx+\frac{1}{p}\int_\Omega
\nabla |x|^{\alpha-2p+2}\cdot\nabla |u|^p~\!dx.
\end{eqnarray*}
Some care is needed in case $p\in(1,2)$.
We first prove (\ref{eq:Rellich_Mit_Omega}) for a fixed function
$u\in C^2_N(\overline\Omega\setminus\{0\})$. For any $\eps>0$,
we define
\begin{gather*}
\f_{\eps}=\left(|u|^2+\eps^2\right)^{\frac{p-2}{2}}u~,\quad \hat\f_\eps=|x|^{\alpha-2p+2}~\!\f_\eps\\
\Phi_\eps=\left(|u|^2+\eps^2\right)^{\frac{p}{4}}-\eps^{p/2}~,\quad
\Theta_\eps=\left(|u|^2+\eps^2\right)^{\frac{p}{2}}-\eps^{p}.
\end{gather*}
Notice that
$$
\nabla \f_\eps\cdot\nabla u\ge(p-1)\left(|u|^2+\eps^2\right)^{\frac{p-4}{2}}|u|^2|\nabla u|^2\ge \frac{4(p-1)}{p^2}|\nabla \Phi_\eps|^2.
$$
In addition, the Hardy inequality gives
\begin{eqnarray*}
\int_\Omega|x|^{\alpha-2p+2}|\nabla \Phi_\eps|^2~\!dx&\ge& \left(\frac{n-2p+\alpha}{2}\right)^2\int_\Omega|x|^{\alpha-2p}|\Phi_\eps|^2~\!dx\\
&=&\left(\frac{n-2p+\alpha}{2}\right)^2\int_\Omega|x|^{\alpha-2p}|u|^p~\!dx+o(1)
\end{eqnarray*}
as $\eps\to 0$. Therefore
$$
\int_\Omega|x|^{\alpha-2p+2}\nabla \f_\eps\cdot\nabla u~\!dx\ge
(p-1)H_{2,\alpha}^2 
\int_\Omega|x|^{\alpha-2p}|u|^p~\!dx+o(1).
$$
We notice that $p~\!\f_\eps\nabla u=\nabla\Theta_\eps$ and we integrate by parts to compute
\begin{eqnarray*}
\int_\Omega\nabla |x|^{\alpha-2p+2}\cdot(\f_\eps\nabla u)~\!dx&=&\!\!-\frac{1}{p}
\int_\Omega(\Delta |x|^{\alpha-2p+2})\Theta_\eps~\!dx\\
&=&\!\!-\frac{(\alpha-2p+2)(n-2p+\alpha)}{p}\int_\Omega|x|^{\alpha-2p}\Theta_\eps~\!dx\\
&=&\!\!-{(\alpha-2p+2)}H_{2,\alpha}\int_\Omega|x|^{\alpha-2p}|u|^p~\!dx+o(1).
\end{eqnarray*}
Now we use integration by parts and H\"older inequality to estimate
\begin{eqnarray*}
\int_\Omega\nabla u\cdot\nabla(|x|^{\alpha-2p+2}\f_\eps)~\!dx&=&\!\!\!
\int_\Omega (-\Delta u)(|x|^{\alpha-2p+2}\f_\eps)~\!dx\\
&\le&\!\!\!\left(\int_\Omega|x|^\alpha|\Delta u|^p~\!dx\right)^{\frac{1}{p}}
\left(\int_\Omega|x|^{\alpha-2p}|\f_\eps|^\frac{p}{p-1}~\!dx\right)^{\frac{p-1}{p}}\\
&\le&\!\!\!\left(\int_\Omega|x|^\alpha|\Delta u|^p~\!dx\right)^{\frac{1}{p}}\!\!
\left(\int_\Omega|x|^{\alpha-2p}|u|^p~\!dx\right)^{\frac{p-1}{p}}\!\!\!+o(1).
\end{eqnarray*}
Since
$$
\int_\Omega\nabla u\cdot\nabla( |x|^{\alpha-2p+2}\f_\eps)~\!dx
=\int_\Omega|x|^{\alpha-2p+2}\nabla \f_\eps\cdot\nabla u~\!dx+
\int_\Omega
\nabla |x|^{\alpha-2p+2}\cdot(\f_\eps\nabla u)~\!dx,
$$
by gluing all above information  we infer
$$
\left(\int_\Omega|x|^\alpha|\Delta u|^p~\!dx\right)^{\frac{1}{p}}\!\!
\left(\int_\Omega|x|^{\alpha-2p}|u|^p~\!dx\right)^{\!\!\frac{p-1}{p}}
\ge \gamma_{p,\alpha}\int_\Omega|x|^{\alpha-2p}|u|^p~\!dx+o(1)~\!,
$$
and (\ref{eq:Rellich_Mit_Omega}) readily follows for $u\in C^2_N(\overline\Omega\setminus\{0\})$, as $\gamma_{p,\alpha}>0$.
To prove (\ref{eq:Rellich_Mit_Omega}) for $C^2_N(\overline\Omega)$ notice that 
$|x|^\alpha,|x|^{\alpha-2p}\in L^1_{\rm loc}(\Omega)$
and use an approximation argument.

It remains to check that the constant $\gamma_{p,\alpha}$ can not be improved if
$0\in\Omega$. Take a nontrivial  function $\omega\in C^2_c(0,1)$, and then
use $u_\eps(x)=|x|^{-H}\omega(|x|^\eps)$ as test function, where 
$H=H_{2,\alpha}$ and $\eps>0$ is a small parameter,
so that $u_\eps\in C^2_c(\Omega\setminus\{0\})$. Then compute
$$
\frac{\displaystyle
\int_\Omega |x|^{\alpha}|\Delta u_\eps|^p dx}
{\displaystyle\int_\Omega |x|^{\alpha-2p}|u_\eps|^{p }dx}=
\frac
{\displaystyle\int_0^1 s^{-1}|\eps^2 s^2\omega''+\eps s(n-2-2H+\eps)\omega'-
\gamma_{p,\alpha}\omega|^p~\!ds}
{\displaystyle\int_0^1 s^{-1}|\omega|^p~\!ds},
$$
let $\eps\to 0$ and conclude.
\QED
Now let $\alpha\in\R$ and assume $\mu_{p,\alpha}>0$. Define the space $\mathcal D^{2,p}(\R^n;|x|^\alpha dx)$ as the completion of $C^{2}_{c}(\R^n\setminus\{0\})$ 
with respect to the norm
\[\|u\|^p_{2,\alpha} =\int_{\R^n} |x|^\alpha |\Delta u|^p dx.\]
Then $\mathcal D^{2,p}(\R^n;|x|^\alpha dx)$ is continuously embedded into $L^{p}(\R^n;|x|^{\alpha-2p} dx)$ and 
$$
 \mu_{p,\alpha} =  \inf_{\scriptstyle u\in \mathcal D^{2,p}(\R^n;|x|^\alpha dx)\atop\scriptstyle  u\ne 0}\frac{\displaystyle
\irn |x|^{\alpha}|\Delta u|^p dx}
{\displaystyle\irn |x|^{\alpha-2p}|u|^{p }dx}~\!.
$$
We introduce the transform
$$
\mathcal T_{2,\alpha}:C^1_c( \mathcal Z^n)\to C^1_c(\R^n\setminus\{0\})~,\quad
(\mathcal T_{2,\alpha}g)(x)=|x|^{-H_{2,\alpha}}g\left(-\log|x|,\frac{x}{|x|}\right).
$$
The "radial version" of the next lemma has been crucially used in \cite{M}.

\begin{Lemma}
\label{L:description2}
Assume that $-\gamma_{p,\alpha}$ is not an eigenvalue
of the Laplace-Beltrami operator
on the sphere. Then
$$
\|g\|_{2,\alpha}:=\left(\irn|x|^\alpha|\Delta(\mathcal T_{2,\alpha} g)|^p~dx\right)^{1/p}
~\quad\text{for $g\in C^2_c( \mathcal Z^n)$,}
$$
is equivalent to the standard norm in $W^{2,p}(\mathcal Z^n)$. Thus $\mathcal T_{2,\alpha}$
can be uniquely extended to an isomorphism   $W^{2,p}(\mathcal Z^n)\to \mathcal D^{2,p}(\R^n;|x|^\alpha dx)$, and
$$
\mathcal D^{2,p}(\R^n;|x|^\alpha dx)=\left\{u\in L^{p}(\R^n;|x|^{\alpha-2p} dx)~|~-\Delta u\in L^p(\R^n;|x|^{\alpha}dx)~\right\}.
$$
If in addition $a>2p-n$, then $C^2_c(\R^n)\subset \mathcal D^{2,p}(\R^n;|x|^\alpha dx)$. 
\end{Lemma}

\Proof
It turns out that $\mu_{p,\alpha}> 0$ by (\ref{eq:MSS}).
By direct computation one has that
\begin{equation}
\label{eq:EF_role}
\begin{cases}
~\displaystyle\irn |x|^{\alpha}|\Delta u|^{p} dx=\displaystyle\int_{\mathcal Z^n}
\left|\Delta_\sigma g+g_{ss}-2A_{p,\alpha}g_{s}-\gamma_{p,\alpha} g\right|^{p}~dsd\sigma\\
~\\
~\displaystyle\irn|x|^{\alpha-2{p}}| u|^{p} dx=\displaystyle\int_{\mathcal Z^n}
|g|^{p}~dsd\sigma~\!,
\end{cases}
\end{equation}
where 
$$
A_{p,\alpha}=\frac{n+2}{2}-\frac{n+\alpha}{p}.
$$
To conclude, adapt the arguments in the proof of Lemma \ref{L:description1}.
\QED

\begin{Remark}
\label{R:mu_symmetry}
The function $\alpha\mapsto\mu_{p,\alpha}$ is symmetric with respect to
$$
\alpha^*_p:=p+n~\!\frac{p-2}{2}.
$$
In fact,  for any $u\in C^2_c(\R^n\setminus\{0\})$, $t\in\R$ it turns out that
\begin{gather*}
\irn|x|^{\alpha^*_p-t}|\Delta \hat u|^p~\!dx=\irn|x|^{\alpha^*_p+t}|\Delta  u|^p~\!dx\\
\irn|x|^{\alpha^*_p-t-2p}| \hat u|^p~\!dx=\irn|x|^{\alpha^*_p+t-2p}|  u|^p~\!dx,
\end{gather*}
where
$$
\hat u(x):=|x|^{\frac{2t}{p}}~\!u\left(\frac{x}{|x|^2}\right).
$$
The weighted second order Emden-Fowler transforms can be used to avoid
boring computations.  Indeed, putting 
$\hat g:= \mathcal T^{-1}_{2,\alpha^*_p-t}\hat u$, $g=\mathcal T^{-1}_{2,\alpha^*_p+t} u$, one has that
$\hat g(s,\sigma)= g(-s,\sigma)$.
To conclude, use (\ref{eq:EF_role}) and notice that the function 
$\alpha\mapsto \gamma_{p,\alpha}$ is even with respect to $\alpha^*_p$
while the function $\alpha\mapsto A_{p,\alpha}$ is odd with respect to $\alpha^*_p$,
that is,
$\gamma_{p,\alpha^*_p-t}=\gamma_{p,\alpha^*_p+t}$ and $A_{p,\alpha^*_p-t}=-A_{p,\alpha^*_p+t}$
for any $t\in\R$.

In particular, we have that $\mu_{p,\alpha}\neq 0$ if and only if $\mu_{p,\hat\alpha}\neq 0$,
where $\hat\alpha=2\alpha^*_p-\alpha$.
In this case, the spaces $\mathcal D^{2,p}(\R^n;|x|^\alpha dx)$ and 
$\mathcal D^{2,p}(\R^n;|x|^{\hat \alpha} dx)$ can be identified trough the isometry
$u\mapsto \hat u$.
\end{Remark}

The next corollary is an immediate consequence of 
Lemma \ref{L:description2} and of  Sobolev embedding theorems for
$W^{2,p}(\mathcal Z^n)$. 

\begin{Corollary}
\label{C:deltanabla}
Let  $p,q$ be given exponents, such that 
$1<p\le q<\infty$ and $q\le \pstarstar$ if $n>2p$. Let $\alpha\in\R$ and assume that
$-\gamma_{p,\alpha}$ is not an eigenvalue
of the Laplace-Beltrami operator
on the sphere. Then

$~i)$  $\mathcal D^{2,p}(\R^n;|x|^\alpha dx)$
is continuously embedded into
$ L^q(\R^n;|x|^{-n+q\frac{n-2p+\alpha}{p}} dx)$.

$ii)$ 
 $\mathcal D^{2,p}(\R^n;|x|^\alpha dx)$
is continuously embedded into
$\mathcal D^{1,p}(\R^n;|x|^{{-n+q\frac{n-p+\alpha}{p}}}dx)$.
\end{Corollary}

Corollary \ref{C:deltanabla} readily implies
that $S_{p,q}(\alpha)>0$, compare with (\ref{eq:space_weight}). Notice that
the function $\alpha\mapsto S_{p,q}(\alpha)$ 
is symmetric with respect to
$\alpha^*_p=p+n~\!\frac{p-2}{2}$ by Remark \ref{R:mu_symmetry}.

\begin{Remark}
\label{R:density_0}
Assume $n>2p$ and take $\alpha=0$. From Lemma \ref{L:description2} we infer that
$ \mathcal D^{2,p}(\R^n;|x|^0 dx)$ coincides with the standard space $ \mathcal D^{2,p}(\R^n)$.
In particular, $ \mathcal D^{2,p}(\R^n)$ can be identified with the standard Sobolev space
$W^{2,p}(\R^n\times \mathbb S^{n-1})$ trough the transform $\mathcal T_{2,\alpha}$, and
$$
\mathcal D^{2,p}(\R^n)=\left\{u\in L^p(\R^n;|x|^{-2p}dx)~|~-\Delta u\in L^p(\R^n)~\right\}.
$$
\end{Remark}

\subsection{Bounded domains}
\label{SS:bounded}
~\\
Here we assume that  $\Omega\subset\R^n$ is a bounded domain of class $C^2$
containing the origin.
We start with a density lemma.

\begin{Lemma}
\label{L:Omega_L}
Assume $\alpha>2p-n$ and let
$u\in C^2_N(\overline\Omega)$. Then there exists
a sequence $u_h\in C^2_N(\overline\Omega\setminus\{0\})$
such that $\Delta u_h\to \Delta u$ in $L^p(\Omega;|x|^\alpha dx)$
and $u_h\to  u$ in $L^q(\Omega;|x|^{-\beta} dx)$, for any
$q\ge p$ and $\beta$ as in (\ref{eq:ma_a}).
\end{Lemma}

\Proof
Take a smooth function $\eta:\R\to \R$ such that
$0\le \eta\le 1$, $\eta(s)\equiv 1$ for $s\le 1$ and $\eta\equiv 0$ for $s\ge 2$. 
We put  $\eta_h(x)=\eta\left(-h^{-1}{\log|x|}\right)$ and we check by direct computation that the sequence
$u_{h}=\eta_hu$ satisfies the desired requirements. Notice that
$\eta_h-1, \nabla \eta_h,\Delta \eta_h\to 0$ pointwise on $\R^n\setminus\{0\}$ and
that the sequences 
$\eta_h, |x||\nabla \eta_h|$ and $|x|^2|\Delta\eta_h|$ are uniformly bounded. In addition,
the supports of $\eta_h-1, \nabla\eta_h$ and $ \Delta\eta_h$ are contained in the  
closed ball of radius $e^{-h}$
about the origin. Thus
\begin{gather*}
\eta_h\to 1~~\text{in $L^q(\Omega;|x|^\nu dx)$ for any $q\ge 1$ and $\nu>-n$,}\\
|\nabla \eta_h|\to 0~~\text{in $L^p(\Omega;|x|^\nu dx)$ for any $\nu>p-n$}~,\quad
\Delta\eta_h\to 0~~\text{in $L^p(\Omega;|x|^\alpha dx)$.}
\end{gather*}
In particular, $u_h\to u$ in $L^q(\Omega;|x|^\nu dx)$ for any $q\ge 1$, $\nu>-n$, and since
$$
\Delta u_h=\eta_h\Delta u+2\nabla \eta_h\cdot\nabla u+u\Delta\eta_h~\!,
$$
then $\Delta u_h\to \Delta u$ in $L^p(\Omega;|x|^\alpha dx)$.
\QED

Now assume that $\alpha$ satisfies
(\ref{eq:ma_b}), and  define
$\mathcal D^{2,p}_{N}(\Omega;|x|^\alpha dx)$ as the completion of $C^2_N(\overline\Omega)$ with respect to the norm
$$
\|u\|_{2,\alpha}=\left(\int_\Omega|x|^\alpha |\Delta u|^p dx\right)^{1/p}.
$$
Then $\mathcal D^{2,p}_N(\Omega;|x|^\alpha dx)$ is continuously embedded into
$L^p(\Omega;|x|^{\alpha-2p})$ by Lemma \ref{L:Mit}.

 Since $\partial\Omega$ is smooth and compactly
contained in $\R^n\setminus\{0\}$, and since 
$C^2_N(\overline\Omega\setminus\{0\})$ is dense in $\mathcal D^{2,p}_N(\Omega;|x|^\alpha dx)$
by Lemma \ref{L:Omega_L}, then $\mathcal D^{2,p}_N(\Omega;|x|^\alpha dx)$ is the space of functions
$u\in L^p(\Omega;|x|^{\alpha-2p})$ such that $\Delta u\in L^{p}(\Omega;|x|^\alpha dx)$
and $u= 0$ in the sense of traces on $\partial\Omega$, coherently with the
definition already given in the introduction.

\section{A linear problem on $\R^n$}
\label{S:linear}
Here we deal with the non-homogeneous equation
\begin{equation}
\label{eq:linear}
-\Delta v=f\quad\text{in $\R^n$.}
\end{equation}
In this section we will always assume that
(\ref{eq:ma_b}) is satisfied.
We start with a simple lemma and then we prove an existence result.

\begin{Lemma}
\label{L:unique}
If $v\in L^p(\R^n;|x|^{\alpha-2p} dx)$ is harmonic on $\R^n\setminus\{0\}$, then
$v\equiv 0$.
\end{Lemma}

\Proof
First of all we notice that $v\in L^1_{\rm loc}(\R^n)$, argue as in Remark \ref{R:L1}.
Next, fix any
$\eta\in C^\infty_c(\R^n)$ and put
$$
p'=\frac{p}{p-1}~,\quad \widetilde\alpha=\frac{2p-\alpha}{p-1}=2p'-\frac{\alpha}{p-1}.
$$ 
Notice that $2p'-n<\widetilde\alpha<np'-n$, as $2p-n<\alpha<np-n$. Since the weights
$|x|^{\widetilde\alpha}$ and $|x|^{\widetilde\alpha-p'}$ are in $L^1_{\rm loc}(\R^n)$, then 
clearly $\eta\in \mathcal D^{2,p'}(\R^n;|x|^{\widetilde\alpha}dx)$, compare with Lemma \ref{L:description2}. Thus
there exists a sequence $\eta_h\in C^2_c(\R^n\setminus\{0\})$ such that
$\eta_h\to \eta$  in ~$\mathcal D^{2,p'}(\R^n;|x|^{\widetilde\alpha} dx)$. 
Since
$$
0=\irn v\Delta \eta_h~\!dx=\irn \left(|x|^{\frac{\alpha-2p}{p}}v\right)
\left(|x|^{\widetilde\alpha p'}\Delta\eta_h\right)~\!dx=\irn v\Delta\eta~\!dx+o(1),
$$
we readily infer that $-\Delta v=0$ on $\R^n$, as $\eta$ was arbitrarily chosen. Thus
$v$ is the null function in $\mathcal D^{2,p}(\R^n;|x|^\alpha dx)$ by Lemma \ref{L:description2}.
\QED

\begin{Theorem}
\label{T:linear}
Let $p>1$, $\alpha\in(2p-n,np-n)$ and let $f\in L^p(\R^n;|x|^\alpha dx)$ be a 
given function.
Then there exists a unique $v\in  \mathcal D^{2,p}(\R^{n};|x|^\alpha dx)$ that solves
(\ref{eq:linear}) in the distributional sense on $\R^n$. If in addition $f\neq 0$ and $f\ge 0$
almost everywhere on $\R^n$, then $v$ is superharmonic and strictly positive on $\R^n$.
\end{Theorem}

\Proof
For any $R>1$
we denote by $A_R$ the annulus $\{R^{-1}<|x|<R\}$. Notice that
$f\in L^p(A_R)$. Let 
$h>1$ be an integer and  let $v_h\in W^{2,p}_N(A_h)$ be the unique solution to
\begin{equation}
\label{eq:vh_ext2}
\begin{cases}
-\Delta v_{h}= f&\text{in $A_h$}\\
v_h=0&\text{on $\partial A_h$.}
\end{cases}
\end{equation}
Now we extend $v_{h}$ by the null function outside $A_h$ and
we write $v_h$ instead of $\chi_{A_h} v_h$ to simplify notation. We use
the Rellich inequality (\ref{eq:Rellich_Mit_Omega}) to infer 
$$
\int_{\R^{n}}|x|^{\alpha-2p}|v_{h}|^{p}dx=\int_{A_h}|x|^{\alpha-2p}|v_{h}|^{p}dx\le
c \int_{A_h}|x|^\alpha|\Delta v_{h}|^{p}dx\le c \irn|x|^\alpha|f|^{p}dx,
$$
where $c=\gamma_{p,\alpha}^{-p}>0$. Therefore, we have that  the sequence $(v_{h})$ is uniformly bounded in $L^{p}(\R^{n},|x|^{\alpha-2p}dx)$, and we can assume
that $v_h\weak v$ weakly in $L^{p}(\R^{n},|x|^{\alpha-2p}dx)$ for some
$v\in L^{p}(\R^{n},|x|^{\alpha-2p}dx)$. For every fixed $R>1$
and for $h>R$ we clearly have
$$
\int_{A_{R}}|\Delta v_{h}|^{p}~\!dx\le c_R\irn |x|^\alpha |f|^{p}~\!dx,
$$
where $c_R$ denotes any constant that might depend on $R$ and $\alpha$ but not on $h$.  
Since $(v_h)$ is  bounded in $L^p(A_R)$ then 
$(v_{h})$ is bounded in $W^{2,p}(A_{R})$. In particular, from every subsequence we can extract a new subsequence 
$v_{h_j}$ such that $v_{h_j}\weak v_{R}$ weakly in $W^{2,p}(A_{R})$ for some $v_{R}\in W^{2,p}(A_{R})$. 
Actually, since $v_{h_j}\to v_{R}$ strongly in $L^{p}(A_{R})$, we have $v_{R}=v$ a.e. on $A_{R}$.
Thus $v\in W^{2,p}_{\rm loc}(\R^n\setminus\{0\})$, $v_h\weak v$ in
$W^{2,p}_{loc}(\R^{n}\setminus\{0\})$, and $v$ solves (\ref{eq:linear})
almost everywhere on $\R^n$. Since $v\in W^{2,p}_{\rm loc}(\R^n\setminus\{0\})$, 
then $v$ has a distributional Laplacian $-\Delta v\in \mathcal D'(\R^n\setminus\{0\})$,
and $v$ solves (\ref{eq:linear}) in the distributional 
sense on $\R^n\setminus\{0\}$.
Arguing as in the proof of Lemma \ref{L:unique}, one gets that
$v$ solves (\ref{eq:linear}) in the $\mathcal D'(\R^n)$ sense,
and in particular $-\Delta v\in L^p(\R^n;|x|^\alpha dx)$.
Since in addition
$v\in L^p(\R^n;|x|^{\alpha-2p}~dx)$, from Lemma \ref{L:description2} we infer that
$v\in  \mathcal D^{2,p}_{N}(\Omega;|x|^\alpha dx)$.
The uniqueness of $v$ readily  follows by Lemma \ref{L:unique}.

The last conclusion in case $f\ge 0$ is immediate, as
$v\in L^1_{\rm loc}(\R^n)$.
\QED

As a consequence of the above results we get a new characterization of
the space $ \mathcal D^{2,p}(\R^n;|x|^\alpha dx)$.

\begin{Corollary}
\label{C:hole}
If $\omega\in   L^p(\R^{n};|x|^{\alpha-2p} dx)$ solves $-\Delta\omega= f$
in the distributional sense on $\R^n\setminus\{0\}$ for some
$f\in  L^p(\R^{n};|x|^{\alpha} dx)$, then $\omega\in \mathcal D^{2,p}(\R^n;|x|^\alpha dx)$.
\end{Corollary} 

\Proof
Use  Theorem \ref{T:linear} to find $v\in  \mathcal D^{2,p}(\R^n;|x|^\alpha dx)$ 
such that $-\Delta v=f$ on $\R^n$. Then $v-\omega\in 
L^p(\R^n;|x|^{\alpha-2p} dx)$ is harmonic on $\R^n\setminus\{0\}$.
Hence, $\omega=v$ by Lemma \ref{L:unique}.
\QED

\section{Proof of theorem \ref{T:noBP}}
\label{SS:weights}

Use Theorem \ref{T:linear} to find the unique superharmonic and positive function $v\in \mathcal D^{2,p}(\R^n;|x|^\alpha dx)$ such that
$$
-\Delta v=|\Delta u|\quad\text{on $\R^n$.}
$$ 
Since
$v\pm u\in \mathcal D^{2,p}(\R^n;|x|^\alpha dx)$, then $v\pm u\in\mathcal 
D^{1,p}(\R^n;|x|^{\alpha-p}dx)$ 
by Corollary \ref{C:deltanabla}. In addition $-\Delta(v\pm u)\ge 0$ on $\R^n$, that implies
$v\pm u\ge 0$ on $\R^n$ by Theorem \ref{T:MM}. Thus $v\ge |u|$ a.e. on $\R^n$, and therefore
\begin{eqnarray*}
S_{p,q}(\alpha)\left(\displaystyle\irn |x|^{-\beta}v^q~\!dx\right)^{p/q}\!\!\!&\le&
\displaystyle\irn|x|^\alpha|\Delta v|^p~\!dx\\
&=&
\displaystyle\irn|x|^\alpha|\Delta u|^p~\!dx=
S_{p,q}(\alpha)\left(\displaystyle\irn |x|^{-\beta}|u|^q~\!dx\right)^{p/q}\\
&\le&S_{p,q}(\alpha)\left(\displaystyle\irn |x|^{-\beta}v^q~\!dx\right)^{p/q},
\end{eqnarray*}
that readily gives $|u|=v$, as $v^q-|u|^q\ge 0$. Since $v>0$ then
up to a change of sign we can assume that $u=v$. In particular, $u$ is superharmonic
and positive.
\QED

\section{Proof of Theorem \ref{T:RS_Omega}}
\label{SS:alpha_Navier}

Fix any $u\in C^2_{N}(\overline\Omega\setminus\{0\})\setminus\{0\}$ and use
Theorem \ref{T:linear} to find $v\in \mathcal D^{2,p}(\R^n;|x|^\alpha dx)$, $v>0$, such that
$$
-\Delta v=\chi_\Omega|\Delta u|\quad\text{on $\R^n$.}
$$
Since $v\pm u\in L^1(\Omega)$, $-\Delta(v\pm u)\ge 0$ in $\Omega$ and
$v\pm u\ge 0$ on $\partial\Omega$, then 
$v\ge |u|$ on $\Omega$. In particular
$$
S_{p,q}(\alpha)\le \frac{\displaystyle\irn|x|^{\alpha}|\Delta v|^p~\!dx}
{\left(\displaystyle\irn |x|^{-\beta}|v|^q~\!dx\right)^{p/q}}\le 
\frac{\displaystyle\int_\Omega|x|^{\alpha}|\Delta u|^p~\!dx}
{\left(\displaystyle\int_\Omega |x|^{-\beta}|u|^q~\!dx\right)^{p/q}}~\!.
$$
Thus $S_{p,q}(\alpha)\le S_{p,q}(\Omega;\alpha)$, as $u$ was arbitrarily chosen
and thanks to the result in Subsection \ref{SS:bounded}. 
Next notice that
$$
S_{p,q}(\Omega;\alpha)\le
\inf_{\scriptstyle u\in C^{2}_{c}(\Omega)\atop\scriptstyle  u\ne 0}
\frac{\displaystyle\int_\Omega |x|^{\alpha}|\Delta u|^p dx}
{\left(\displaystyle\int_\Omega |x|^{-\beta}|u|^{q}~\!dx\right)^{p/q}}=S_{p,q}(\alpha)
$$
by simple inclusion and rescaling
arguments, and (\ref{eq:weight_Navier}) is proved. 

It remains to prove that 
$S_{p,q}(\Omega;\alpha)$ is not achieved.
Assume by contradiction that
there exists $u\neq 0$ in $\mathcal D^{2,p}_{N}(\Omega;|x|^\alpha dx)$ achieving  $S_p^{\rm Nav}(\Omega)$, and define $v\in \mathcal D^{2,p}(\R^n;|x|^\alpha dx)$ as before. Then $v>0$,
$v\ge \chi_\Omega|u|$ and from $-\Delta v= \chi_\Omega|\Delta u|$ we infer
\begin{eqnarray*}
S_p\left(\displaystyle\irn |x|^{-\beta}v^q~\!dx\right)^{p/q}&\le&
\displaystyle\irn|x|^\alpha|\Delta v|^p~\!dx\\
&=&
\displaystyle\int_\Omega|x|^\alpha|\Delta u|^p~\!dx=
S_p\left(\displaystyle\int_\Omega |x|^{-\beta}|u|^q~\!dx\right)^{p/q}
\end{eqnarray*}
by (\ref{eq:equality}). Thus
$$
\irn |x|^{-\beta}\left(v^q-|\chi_\Omega u|^q\right)~\!dx\le 0,
$$
that together with $v^q-|\chi_\Omega u|^q\ge 0$ implies $|\chi_\Omega u|=v$. Clearly this is impossible, as
$v>0$ in $\R^n$.
\QED

\section{Proof of Theorem \ref{T:Rellich}}
\label{S:Rellich}
Up to a dilation and thanks to Sobolev embedding theorem, we can assume that 
$\overline B_1\subset\Omega$ and $q\le \pstarstar$ if $n>2p$.
One can adapt the choice of  test functions that has been made in
Remark \ref{R:resonance}
to check that $S_{p,q}^{\rm Nav}(\Omega;\alpha)=0$
if $\alpha=np-n$.  Thus it suffices to prove the result in case $\alpha>np-n$.
 
For $X=W^{2,p}_N(\Omega)$, $X=C^2_N(\overline\Omega)$ or $X= C^2_N(\overline\Omega\setminus\{0\})$ and for $v\in X$ we put
$$
m(X)=\inf_{\scriptstyle v\in X\atop\scriptstyle  v\ne 0} R(v)~,\quad
R(v):=\frac{\displaystyle
\int_\Omega |x|^{\alpha}|\Delta v|^p dx}
{\left(\displaystyle\int_\Omega |x|^{-\beta}|v|^{q}dx\right)^{p/q}}~\!.
$$
Since $\alpha>np-n>0$, then $\Delta v\in L^p(\Omega;|x|^\alpha dx)$ for any 
$v\in W^{2,p}(\Omega)$. In particular,
the infima $m(X)$ are well defined. By trivial inclusions and thanks to 
 Lemma \ref{L:Omega_L}  we have that
$$
m(W^{2,p}_N(\Omega))\le m(C^2_N(\overline\Omega))\le
m(C^2_N(\overline\Omega\setminus\{0\}))=S_{p,q}^{\rm Nav}(\Omega;\alpha).
$$
Actually $m(C^2_N(\overline\Omega))=m(C^2_N(\overline\Omega\setminus\{0\}))$
by Lemma \ref{L:Omega_L}. Therefore, to conclude the proof we have to show that
\begin{equation}
m(C^2_N(\overline\Omega))=0. 
\label{eq:claim2}
\end{equation}
The main step consists in proving that
\begin{equation}
\label{eq:claim1}
m(W^{2,p}_N(\Omega))=0.
\end{equation}
Since $\alpha>np-n$, we have that
$\gamma_\alpha:=\gamma_{p,\alpha}<0$; compare with (\ref{eq:gamma_intro}). In particular, we can find a geodesic ball
$\mathcal B\subset \mathbb S^{n-1}$ such that 
$$
-\gamma_{\alpha}=
\inf_{\scriptstyle \f\in H^1_0(\mathcal B)\atop\scriptstyle  \f\ne 0}
\frac{\displaystyle\int_{\mathcal B}|\nabla_\sigma \f|^2~\!d\sigma}
{\displaystyle\int_{\mathcal B}|\f|^2~\!d\sigma}.
$$
Fix an eigenfunction $\f\in H^1_0(\mathcal B)$ relative to the eigenvalue
$-\gamma_{\alpha}$
and any nontrivial function $\omega\in C^2_c(\R_+)$. For any small $\eps>0$ use polar coordinates
$(r,\sigma)\in\R_+\times\mathbb S^{n-1}$ to define 
$$
u_\eps(r\sigma):=r^{-H}\omega(r^\eps)\f(\sigma)~\!,
$$ 
where $H=H_{2,\alpha}$ is defined in (\ref{eq:H}). Let $\Omega_\eps$ be the support of
$u_\eps$. Then $\Omega_\eps$ has a Lipschitz boundary and it is
compactly contained in $\Omega$ for any $\eps$
small enough. In addition, $u_\eps\in W^{1,p}_0\cap W^{2,p}(\Omega_\eps)$. Let 
$v_\eps\in W^{2,p}_N(\Omega)$ be the solution of
$$
\begin{cases}
-\Delta v=\chi_{\Omega_\eps}|-\Delta u_\eps|&\text{in $\Omega$,}\\
v=0&\text{on $\partial{\Omega}$.}
\end{cases}
$$
By Lemma \ref{L:extension1} we have that
$v_\eps\ge|u_\eps|$ in $\Omega_\eps$, and thus
$$
m(W^{2,p}_N(\Omega))\le R(v_\eps)
\le \frac{\displaystyle
\int_{\Omega_\eps} |x|^{\alpha}|\Delta u_\eps|^p dx}
{\left(\displaystyle\int_{\Omega_\eps} |x|^{-\beta}|u_\eps|^{q}dx\right)^{p/q}}.
$$
Since $\Delta_\sigma\f=\gamma_\alpha\f$ on $\mathcal B$, we compute
$$
(\Delta u_\eps)(r\sigma)
=\left[\Delta_r\left(r^{-H}\omega(r^\eps)\right)+\gamma_\alpha r^{-H-2}\omega(r^\eps)
\right]~\!\f(\sigma)
$$
where $r>0$, $\sigma\in \mathcal B$ and $\Delta_r w= w''+(n-1)r^{-1}w'$ for any $w\in C^2_c(\R_+)$. Therefore
$$
(\Delta u_\eps)(x)=\eps r^{-H-2+\eps}\left[
\eps r^\eps\omega''(r^\eps)+(c+\eps)\omega'(r^\eps)\right]~\!\f(\sigma)
$$
where $c$ denotes any constant independent on $\eps$, and
\begin{gather*}
\int_{\Omega_\eps}|x|^\alpha|\Delta u_\eps|^p~\!dx=
c\eps^{p-1}\int_0^\infty s^{p-1}|\eps s\omega''+(c+\eps)\omega'|^p~\!ds\\
\displaystyle\int_\Omega |x|^{-\beta}|u_\eps|^{q}dx=c\eps^{-1}
\displaystyle\int_0^\infty s^{-1}|\omega|^{q}dx.
\end{gather*}
In particular, $R(v_\eps)\to 0$ as $\eps\to 0$, and (\ref{eq:claim1}) is proved.

We are in position to prove (\ref{eq:claim2}). Fix any $\delta>0$ and use (\ref{eq:claim1}) to find a nontrivial function $\hat v\in W^{2,p}_N(\Omega)$
such that $R(\hat v)<\delta$. Take a sequence $v_h\in C^2_N(\overline\Omega)$ such that
$v_h\to \hat v$ in $W^{2,p}_N(\Omega)$. Notice that $|x|^\alpha\in L^\infty(\Omega)$ as $\alpha>0$. Therefore
$$
\lim_{h\to\infty}
\int_{\Omega}|x|^\alpha|\Delta v_h|^p~\!dx= \int_{\Omega}|x|^\alpha|\Delta \hat v|^p~\!dx
$$
by Lebesgue's theorem. Then use  Fatou's lemma to get
$$
\lim_{h\to \infty}\int_{\Omega} |x|^{-\beta}|v_{h_j}|^q dx\ge 
\int_{\Omega} |x|^{-\beta}|\hat v|^q dx
$$
up to a subsequence $v_{h_j}$.
Thus we have that
$$
m(C^2_N(\overline\Omega))\le \lim_{j\to\infty}R(v_{h_j})\le R(\hat v)<\delta,
$$
that proves (\ref{eq:claim2}), as $\delta$ was arbitrarily chosen.
The theorem is completely proved.
\QED

We conclude this paper with partial result for lower exponents $\alpha$. 
We take $p=2$, and we put
$$
\gamma_{\alpha}=\gamma_{2,\alpha}=\left(\frac{n-2}{2}\right)^2-\left(\frac{\alpha-2}{2}\right)^2~,
\quad \overline{\gamma}_\alpha=\left(\frac{n-2}{2}\right)^2+\left(\frac{\alpha+2}{2}\right)^2.
$$
If $q=2$ then $\beta=\alpha-4$ and hence
$$
S_{2,2}^{\rm Nav}(\Omega;\alpha):= \inf_{\scriptstyle u\in  C^2_N(\overline\Omega\setminus\{0\})\atop\scriptstyle  u\ne 0}\frac{\displaystyle
\int_\Omega |x|^{\alpha}|\Delta u|^p dx}
{\displaystyle\int_\Omega |x|^{\alpha-4}|u|^{q}dx}.
$$
In \cite{GM}, \cite{CM1} and \cite{MSS} it has been proved that
$$
S_{2,2}(\alpha)=S_{2,2}(\R^n;\alpha)
=\min_{k\in\mathbb N\cup\{0\}}\left|\gamma_\alpha+k(n-2+k)\right|^2.
$$

\begin{Proposition}
\label{P:CM1}
Assume that $\Omega$ is the unit ball in $\R^n$ and that 
$\alpha\le n$. Then
for every $u\in C^{2}_{c}(\overline\Omega\setminus\{0\})$ one has
\begin{equation}
\label{eq:improved}
\int_{\Omega}|x|^{\alpha}|\Delta u|^{2}~\!dx-S_{2,2}(\alpha)\int_{\Omega}|x|^{\alpha-4}|u|^2~\!dx
\ge\frac{\overline\gamma_{\alpha}}{2}\int_{\Omega}|x|^{\alpha-4}|\log|x||^{-2}|u|^2~\!dx~\!.
\end{equation}
In particular, $S_{2,2}^{\rm Nav}(\Omega;\alpha)>0$ for any $\alpha<n$, such that
$-\gamma_\alpha$ is not an eigenvalue of the Laplace-Beltrami operator on the sphere.
\end{Proposition}

In \cite[Theorem 5.1(i)]{CM1} a similar proposition has been stated. However the
assumption $\alpha\le n$ has been neglected there. In view of Theorem 
\ref{T:Rellich}, the assumption $\alpha\le n$ is clearly needed to have (\ref{eq:improved}).
To prove Proposition \ref{P:CM1}, follow the  computations
in \cite{CM1}.


\label{References}

\end{document}